\documentclass[11pt, a4paper]{amsart}
\usepackage[top=35truemm,bottom=35truemm,left=35truemm,right=35truemm]{geometry}
\usepackage{setspace} 
\usepackage{amsmath, amsfonts,amsthm,amssymb,mathrsfs}
\usepackage{enumerate}
\usepackage[pdftex]{graphicx}
\usepackage{ascmac}
\usepackage[arrow,matrix]{xy}
\usepackage{color}
\usepackage{array}
\usepackage{pifont}
\usepackage{longtable}
\usepackage{mathtools}

\theoremstyle{definition}
\newtheorem{Thm}{Theorem}

\newtheorem{theo}{Theorem}
\newtheorem*{theo*}{Theorem}
\newtheorem{defi}{Definition}

\newtheorem{prop}{Proposition}
\newtheorem{cor}{Corollary}
\newtheorem{lem}{Lemma}

\begin{document}

\title[ ]{Growth rates of 3-dimensional hyperbolic Coxeter groups are Perron numbers}
\author{Tomoshige Yukita}
\address{Department of Mathematics, School of Education, Waseda University, Nishi-Waseda 1-6-1, Shinjuku, Tokyo 169-8050, Japan}
\email{yshigetomo@suou.waseda.jp}
\date{}
\thanks{}

\begin{abstract}
In this paper we consider the growth rates of 3-dimensional hyperbolic Coxeter polyhedra some of its dihedral angles are $\frac{\pi}{m}$ for $m\geq{7}$. 
By combining with the classical result by Parry \cite{Pa} and the main result of \cite{Y}, we prove that the growth rates of 3-dimensional hyperbolic Coxeter groups are Perron numbers. 
\end{abstract}

\maketitle


\setstretch{1.1}
\section{Introduction}
Let $\mathbb{H}^{n}$ denote the upper half-space model of hyperbolic n-space and $\overline{\mathbb{H}}^{n}$ its closure in $\mathbb{R}^n\cup{\{\infty\}}$. 
A convex polyhedron $P\subset{\overline{\mathbb{H}}^{n}}$ of finite volume is called a \textit{Coxeter polyhedron} if all of its dihedral angles are of the form $\frac{\pi}{m}$ for an integer $m\geq{2}$ or $m=\infty$, i.e. the intersection of the respective facets is a point on the boundary $\partial{\mathbb{H}^{n}}$.
The set $S$ of reflections with respect to the facets of $P$ generates a discrete group $\Gamma$, called a \textit{hyperbolic Coxeter group},
and the pair $(\Gamma, S)$ is called the \textit{Coxeter system} associated with $P$. 
Then $P$ becomes a fundamental domain for $\Gamma$. 
If $P$ is compact (resp. non-compact), the hyperbolic Coxeter group $\Gamma$ is called \textit{cocompact} (resp. \textit{cofinite}). 
The \textit{growth series} $f_{S}(t)$ of a Coxeter system $(\Gamma, S)$ is the formal power series $\sum_{k=0}^{\infty}a_{k}t^{k}$ where $a_{k}$ is the number of elements of $\Gamma$ whose word length with respect to $S$ is equal to $k$. 
Then $\tau_{\Gamma} :=\limsup_{k \rightarrow \infty} \sqrt[k]{a_k}$ is called the \textit{growth rate} of the Coxeter system $(\Gamma, S)$. 
By means of the Cauchy-Hadamard theorem, $\tau_{\Gamma}$ is equal to the reciprocal of the radius of convergence $R$ of $f_{S}(t)$. 
The growth series and the growth rate of a hyperbolic Coxeter polyhedron $P$ is defined to be the growth series and the growth rate of the Coxeter system $(\Gamma, S)$ associated with $P$, respectively. 
It is known that the growth rate $\tau_{\Gamma}$ of a hyperbolic Coxeter group is a real algebraic integer bigger than 1 \cite{dlH1}.
In the paper \cite{Y}, we showed that the growth rates of 3-dimensional hyperbolic Coxeter polyhedra whose dihedral angles are of the form $\frac{\pi}{m}$ for $m=2,3,4,5,6$ are always Perron numbers. 
Let us recall that a real algebraic integer $\tau>1$ is called a \textit{Perron number} if all its algebraic conjugates are less than $\tau$ in absolute value.
In general, a 3-dimensional hyperbolic Coxeter polyhedron may have some dihedral angles of the form $\frac{\pi}{m}$, $m\geq{7}$, so that in this paper we consider the growth rates of 3-dimensional hyperbolic Coxeter polyhedra some of its dihedral angles are $\frac{\pi}{m}$ for $m\geq{7}$ and prove Theorem A.

\begin{Thm}
The growth rates of non-compact 3-dimensional hyperbolic Coxeter polyhedra with some dihedral angles $\frac{\pi}{m}$ for $m\geq{7}$ are Perron numbers. 
\end{Thm}

The classical result by Parry \cite{Pa}, together with the main result of \cite{Y}  and Theorem A imply Theorem B.

\begin{Thm}
The growth rates of 3-dimensional hyperbolic Coxeter groups are Perron numbers.
\end{Thm}


Theorem B settles the 3-dimensional case of a conjecture by Kellerhals and Perren \cite{KP}. 
Results in the closer vein have recently been published by Kellerhals and Nonaka \cite{KN} and Komori and Yukita \cite{KY}. 
Also, since the set of growth rates of 3-dimensional hyperbolic Coxeter groups is now shown to comprise only Perron numbers by Theorem B, the minimal growth rate becomes a necessary point to be mentioned; 
Kellerhals determined the minimal growth rate among all 3-dimensional cofinite hyperbolic Coxeter groups \cite{Ke}, 
while Kellerhals and Kolpakov found the minimal growth rate in the case of 3-dimensional cocompact hyperbolic Coxeter groups \cite{KK}.


The organization of the present paper is as follows.
In section 2, we review a useful formula from \cite{Y} which allow us to calculate the growth function $f_{S}(t)$.
Then,  we calculate the growth function $f_{S}(t)$ of a non-compact hyperbolic Coxeter polyhedron some of whose dihedral angles are $\frac{\pi}{m}$ for $m\geq{7}$ in section 3. 


\section{Preliminaries}
In this section,  we introduce some notations and review useful identities in \cite{Y} to calculate the growth functions of hyperbolic Coxeter groups.

\begin{defi}{\rm(Coxeter system, Coxeter graph, growth rate)}

{\rm{(i)}} A {\textit {Coxeter system}} $(\Gamma,S)$ consists of a group $\Gamma$ and a set of generators $S\subset{\Gamma}$, $S=\{s_i\}_{i=1}^m$, with relations $(s_is_j)^{m_{ij}}$ for each $i , j$ , where $m_{ii}=1$ and $m_{ij}\geq{2}$ or $m_{ij}=\infty$ for $i\neq{j}$.
We call $\Gamma$ a {\textit {Coxeter group}}.
For any subset $I\subset{S}$, we define $\Gamma_I$ to be the subgroup of $\Gamma$ generated by $\{s_i\}_{i\in{I}}$. 
Then $\Gamma_I$ is called the {\textit {Coxeter subgroup}} of $\Gamma$ generated by $I$.

{\rm{(ii)}} The Coxeter graph of $(\Gamma,S)$ is constructed as follows: \\
Its vertex set is $S$. 
If $m_{ij}\geq{3} \ (s_i\neq{s_j} \in{S})$, we join the pair of vertices by an edge and label such an edge with $m_{ij}$. 
If $m_{ij}=\infty \ (s_i\neq{s_j} \in{S})$, we join the pair of vertices by an bold edge. 

{\rm{(iii)}} The {\textit{growth series}} $f_{S}(t)$ of a Coxeter system $(\Gamma,S)$ is the formal power series $\sum_{k=0}^{\infty}a_kt^k$ where $a_k$ is the number of elements of $\Gamma$ whose word length with respect to $S$ is equal to $k$. Then $\tau =\limsup_{k \rightarrow \infty} \sqrt[k]{a_k}$ is called the {\textit{growth rate}} of $(\Gamma,S)$.

\end{defi}

A Coxeter group $\Gamma$ is \textit{irreducible} if the Coxeter graph of $(\Gamma,S)$ is connected.
In this paper, we are interested in Coxeter groups which act discretely on the hyperbolic n-space.

\begin{defi}{\rm(hyperbolic polyhedron)} 

A subset $P\subset{\overline{\mathbb{H}}^n}$ is called a \textit{hyperbolic polyhedron} if $P$ can be written as the intersection of finitely many closed half spaces: $P=\cap{H_S}$, where $H_S$ is the closed domain of $\mathbb{H}^n$ bounded by a hyperplane $S$.
\end{defi}

Suppose that $S\cap{T}\neq{\emptyset}$  in $\mathbb{H}^{n}$.
Then we define the {\textit {dihedral angle}} between $S$ and $T$ as follows: let us choose a point $x\in{S\cap{T}}$ and consider the outer-normal vectors $e_S, e_T\in{\mathbb{R}^n}$ of $S$ and $T$ with respect to $P$ starting from $x$. 
Then the dihedral angle between $S$ and $T$ is defined as the real number $\theta\in{[0, \pi)}$ satisfying $\cos{\theta}=-(e_S, e_T)$ where $(\cdot , \cdot)$ denote the Euclidean inner product on $\mathbb{R}^n$.

If $S\cap{T}=\emptyset$ in $\mathbb{H}^n$, then $\overline{S}\cap{\overline{T}}\in{\overline{\mathbb{H}}^{n}}$ is a point at the ideal boundary $\partial{\mathbb{H}^n}$ of $\mathbb{H}^{n}$, and we define the dihedral angle between $S$ and $T$ to be equal to zero.

\begin{defi}{\rm(hyperbolic Coxeter polyhedron)}
 
A hyperbolic polyhedron $P\subset{\overline{\mathbb{H}}^{n}}$ of finite volume is called a {\textrm{hyperbolic Coxeter polyhedron}} if all of its dihedral angles have the form $\frac{\pi}{m}$ for an integer $m\geq{2}$ or $m=\infty$, i.e. the intersection of the respective facets is a point on the boundary $\partial{\mathbb{H}^{n}}$.
\end{defi}

Note that a hyperbolic polyhedron $P\subset{\overline{\mathbb{H}}^{n}}$ is of finite volume if and only if $P\cap{\partial{\mathbb{H}^n}}$ consists only of ideal vertices of $P$.
If $P\subset{\overline{\mathbb{H}}^{n}}$ is a hyperbolic Coxeter polyhedron, the set $S$ of reflections with respect to the facets of $P$ generates a discrete group $\Gamma$. We call $\Gamma$ the {\textit{$n$-dimensional hyperbolic Coxeter group}} associated with $P$. Moreover, if $P$ is compact (resp. non-compact), $\Gamma$ is called {\textit{cocompact}} (resp. {\textit{cofinite}}).

We recall Solomon's formula and Steinberg's formula which are very useful for calculating growth series.

\begin{theo}{\rm (Solomon's formula)}\cite{So} \\
The growth series $f_S(t)$ of an irreducible finite Coxeter group $(\Gamma, S)$ can be written as 
$f_S(t)=[m_1+1,m_2+1, \cdots, m_k+1]$ where $[n] =1+t+ \cdots +t^{n-1}, [m,n]=[m][n]$,etc. and $\{m_1, m_2, \cdots, m_k \}$
is the set of exponents of $(\Gamma, S)$.
\end{theo}

The exponents of irreducible finite Coxeter groups are shown in Table 1 (see \cite{Hu} for details).

\begin{table}[h]
\begin{center}
\caption{Exponents}
\begin{tabular}{|c|c|c|}
\hline 
Coxeter group & Exponents & growth series \\ 
\hline
$A_n$ & $1,2,\cdots,n$ & $[2,3,\cdots,n+1]$ \\ 
\hline
$B_n$ & $1,3,\cdots,2n-1$ & $[2,4,\cdots,2n]$ \\ 
\hline
$D_n$ & $1,3,\cdots,2n-3,n-1$ & $[2,4,\cdots,2n-2][n]$ \\ 
\hline
$E_6$ & 1,4,5,7,8,11 & [2,5,6,8,9,12] \\ 
\hline
$E_7$ & 1,5,7,9,11,13,17 & [2,6,8,10,12,14,18] \\ 
\hline
$E_8$ & 1,7,11,13,17,19,23,29 & [2,8,12,14,18,20,24,30] \\ 
\hline
$F_4$ & 1,5,7,11 & [2,6,8,12] \\ 
\hline
$H_3$ & 1,5,9 & [2,6,10] \\ 
\hline
$H_4$ & 1,11,19,29 & [2,12,20,30] \\ 
\hline
$I_2(m)$ & 1,m-1 & [2,m] \\ 
\hline
\end{tabular}
\end{center}
\end{table}

\begin{theo}{\rm (Steinberg's formula)}\cite{St} \\
Let $(\Gamma, S)$ be a Coxeter group.
Let us denote the Coxeter subgroup of $(\Gamma, S)$ generated by the subset $T\subseteq S$ by $(\Gamma_T,T)$, and denote its growth function by $f_T(t)$.
Set $\mathcal{F}=\{T\subseteq S \;:\; \Gamma_T$ is finite $\}$. Then
$$
\frac{1}{f_S(t^{-1})}=\sum _{T \in \mathcal{F}} \frac{(-1)^{|T|}}{f_T(t)}.
$$
\end{theo}

By Theorem 1 and Theorem 2, the growth series of $(\Gamma, S)$ is represented by a rational function $\frac{p(t)}{q(t)} (p,q\in{\mathbb{Z}[t]})$. 
The rational function $\frac{p(t)}{q(t)}$ is called the \textit{growth function} of $(\Gamma, S)$. 
The radius of convergence $R$ of the growth series $f_{S}(t)$ is equal to the real positive root of $q(t)$ which has the smallest absolute value among all its roots. 

Suppose that $P$ is a hyperbolic Coxeter polyhedron in $\overline{\mathbb{H}}^{3}$ and $v$ be a vertex of $P$. 
Let $F_{1},\cdots, F_{n}$ be adjacent facets of $P$ incident to $v$ and $\frac{\pi}{a_{i}}$ be the dihedral angle between $F_{i}$ and $F_{i+1}$. 
By Andreev's theorem \cite{A}, the number of facets of $P$ incident to $v$ is at most 4 and $a_1,\cdots, a_n$ satisfy the following conditions. 
\begin{eqnarray}
a_{1}=a_{2}=a_{3}=a_{4}=2 & \qquad & \text{if $n=4$}. \\
\dfrac{1}{a_{1}}+\dfrac{1}{a_{2}}+\dfrac{1}{a_{3}}\geq{1}  & & \text{if $n=3$} .
\end{eqnarray}


Note that a vertex $v$ of $P$ is an ideal vertex if and only if $a_{1}=a_{2}=a_{3}=a_{4}=2$ or $\dfrac{1}{a_{1}}+\dfrac{1}{a_{2}}+\dfrac{1}{a_{3}}=1$, and we call such a vertex a \textit{cusp}, for short. 
We shall use the following notation and terminology for the rest of the paper: \\
$\bullet$ If a vertex $v$ of $P$ satisfies the identity $(1)$, we call $v$ a \textit{cusp of type $(2,2,2,2)$}. \\
$\bullet$ If a vertex $v$ of $P$ satisfies the inequality $(2)$, we call $v$ a \textit{vertex of type $(a_1, a_2, a_3)$}. \\
$\bullet$ $v_{2,2,2,2}$ denotes the number of cusps of type $(2,2,2,2)$. \\
$\bullet$ $v_{a_1, a_2, a_3}$ denotes the number of vertices of type $(a_1, a_2, a_3)$. \\
$\bullet$ $V, E, F$ denotes the number of vertices, edges and facets of $P$. \\
$\bullet$ If an edge $e$ of $P$ has dihedral angle $\frac{\pi}{m}$,  we call it \textit{$\frac{\pi}{m}$-edge}. \\
$\bullet$ $e_m$ denotes the number of $\frac{\pi}{m}$-edges. \\
$\bullet$ The growth function $f_{S}(t)$ of the Coxeter group $(\Gamma, S)$ associated with $P$ is called the \textit{growth function of $P$}.\\
$\bullet$ The growth rate of the Coxeter group $(\Gamma, S)$ associated with $P$ is called the \textit{growth rate of $P$}.\\

By Lemma 2 in \cite{Y}, the following identities and inequality hold for $P$. 
\begin{eqnarray} 
 & &V-E+F= 2 .\\
 & & V=v_{2,2,2,2}+\sum_{n\geq{2}}^{}v_{2,2,n}+v_{2,3,3}+v_{2,3,4}+v_{2,3,5}+v_{2,3,6}+v_{2,4,4}+v_{3,3,3} .\\
 & & E=\sum_{n\geq2}^{}e_{n} .\\ 
 & &2e_{2} = 4v_{2,2,2,2}+3v_{2,2,2}+2\sum_{n=3}^{\infty}v_{2,2,n}+v_{2,3,3}+v_{2,3,4}+v_{2,3,5}+v_{2,3,6}+v_{2,4,4} .\\
 & &2e_{3} = 3v_{3,3,3}+2v_{2,3,3}+v_{2,2,3}+v_{2,3,4}+v_{2,3,5}+v_{2,3,6} .\\
 & &2e_{4} = 2v_{2,4,4}+v_{2,2,4}+v_{2,3,4} .
   \end{eqnarray}
 \begin{eqnarray}
 & &2e_{5} = v_{2,2,5}+v_{2,3,5} .\\
 & &2e_{6} = v_{2,2,6}+v_{2,3,6} .\\
 & &2e_{n} = v_{2,2,n}  \qquad \text{$n\geq{7}$} .\\
 & & v_{2,2,2,2}+v_{2,3,6}+v_{2,4,4}+v_{3,3,3}\geq{1} .
\end{eqnarray}

We use these identities and inequality to calculate growth functions of cofinite 3-dimensional hyperbolic Coxeter groups.
By using the following proposition, we show that the growth rates of cofinite 3-dimensional hyperbolic Coxeter groups are Perron numbers.

\begin{prop}\rm{(\cite{KU}, Lemma1)}\\
\label{prop:KU}
 Let $g(t)$ be a polynomial of degree $n \geq 2$ having the form
 $$
 g(t)=\sum _{k=1}^n a_k t^k -1, 
 $$
 where  $a_k$ are non-negative integers.
 We assume that the greatest common divisor of $\{k \in \mathbb{N} \; |\; a_k \neq 0 \}$ is $1$.
 Then there exists a real number $r_0$,  $0<r_0<1$ which is the unique zero of $g(t)$ having the smallest absolute value among all zeros of $g(t)$.
\end{prop}

Here we review some details of the proof of Theorem 3 in \cite{Y}. 
First, we calculate the growth function $f_{S}(t)$ of $P$ by means of Steinberg's formula.
Second, by using the identities (3)$\sim$(11) and the inequality (12), we proved that all of the coefficients of the denominator polynomial of  $f_{S}(t)$ except its constant term are non-negative. 
This observation is a key point to prove that the growth rates of non-compact hyperbolic Coxeter polyhedra some of whose dihedral angles are $\frac{\pi}{m}$ for $m\geq{7}$ are always Perron numbers. 
Finally, by applying Proposition 1 to the denominator polynomial of $f_{S}(t)$, we conclude that the growth rate of $P$ is a Perron number. 
 

\section{non-compact Coxeter polyhedra some of whose dihedral angles are $\frac{\pi}{m}$ for $m\geq{7}$}
In this section, we calculate the growth function $f_{S}(t)$ of a non-compact hyperbolic Coxeter polyhedron $P$ some of whose dihedral angles are $\frac{\pi}{m}$ for $m\geq{7}$ and prove the growth rate of $P$ is a Perron number.
\begin{theo}
Let the sum of the numbers of the $\frac{\pi}{m}$-edges for $m\geq{7}$ be $k$, that is, 
\[
k=\displaystyle \sum_{m\geq{7}}^{}e_m
\]
 Then we obtain the following inequality.
\[
k\leq{F-3}
\]
Moreover, if the equality $k=F-3$ holds, then $P$ has a unique cusp of type $(2,2,2,2)$ and all other vertices are of type $(2,2,m_{i})$ $(m_{i}\geq{7})$.
\end{theo}

In order to prove Theorem 3, we use the following deformation argument for Coxeter polyhedra introduced by Kolpakov in \cite{K}. 
We present it in a modified form which is more suitable for further account. 

\begin{theo} (\cite{K}, Proposition 1 and 2.)

(i) Suppose that a hyperbolic Coxeter polyhedron $P$ has some $\frac{\pi}{m}$-edges for $m\geq{7}$.
Then all of the $\frac{\pi}{m}$-edges can be contracted to cusps of type $(2,2,2,2)$.
The hyperbolic Coxeter polyhedron $\hat{P}$ which is obtained from $P$ by contracting some $\frac{\pi}{m}$-edges of $P$ is called the {\textit{pinched Coxeter polyhedron }} of $P$. 

(ii) If a hyperbolic Coxeter polyhedron $P$ has some cusps of type $(2,2,2,2)$, then there exists Coxeter polyhedron  which is obtained by opening cusps of type $(2,2,2,2)$. (see Fig 1)

\end{theo}

\begin{figure}[htbp]
\begin{center}
 \includegraphics [width=200pt, clip]{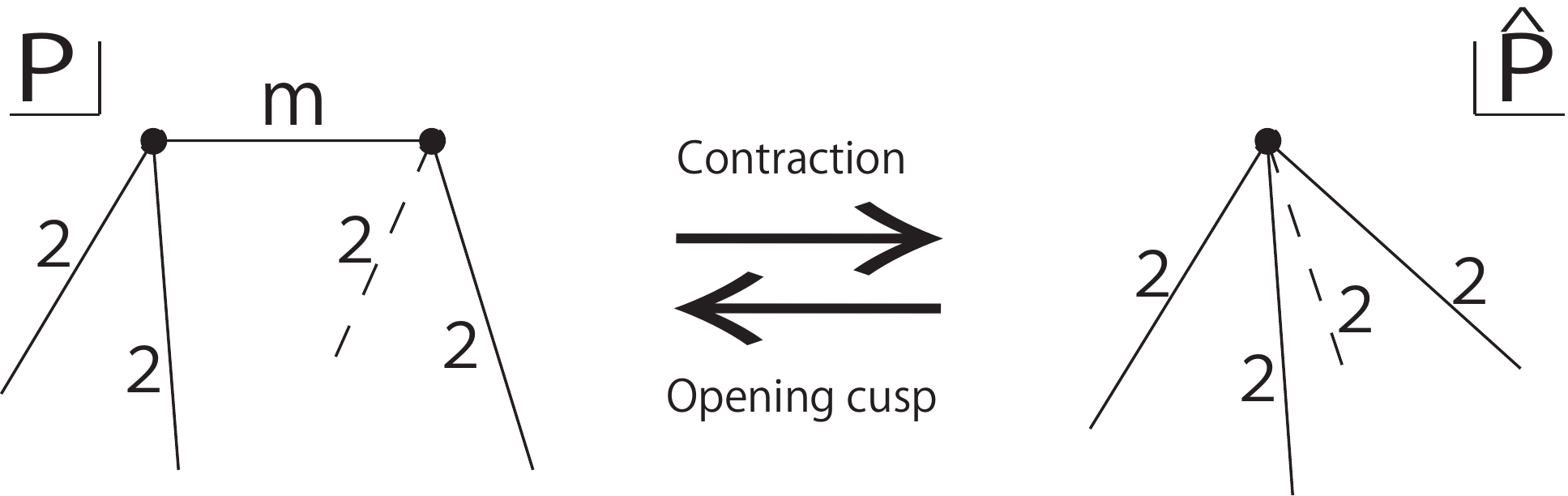}
\end{center}
\caption{}
\label{fig1}
\end{figure} 

In the sequel, $\hat{P}$ denotes the pinched Coxeter polyhedron obtained from $P$ and $\hat{V}, \hat{E}, \hat{F}$, 
$\hat{v}_{2,2,2,2}, \hat{v}_{a_1, a_2, a_3}$ and $\hat{e}_m$ denote respectively the number of vertices, edges, facets, cusps of type $(2,2,2,2)$, vertices of type $(a_1, a_2, a_3)$ and $\frac{\pi}{m}$-edges of $\hat{P}$. 

\textit{Proof of Theorem 3.} Suppose that $P$ is a non-compact hyperbolic Coxeter polyhedron and the sum of the numbers of the $\frac{\pi}{m}$-edges of $P$ is $k$.
The inequalities (3), (4) and (5) imply that 
\begin{eqnarray}
F-2&=&v_{2,2,2,2}+\dfrac{1}{2}\text{(the number of vertices of $P$ with valency 3)} \\ 
\hat{F}-2 &=& \hat{v}_{2,2,2,2}+\dfrac{1}{2}\text{(the number of vertices of $\hat{P}$ with valency 3)} .
\end{eqnarray}
Since $\hat{P}$ is the pinched Coxeter polyhedron of $P$,  $F=\hat{F}$ and $\hat{v}_{2,2,2,2}=v_{2,2,2,2}+k$.
Then the equality (14) implies that 
\begin{eqnarray}
F-2 &=& v_{2,2,2,2}+k+\dfrac{1}{2}\text{(the number of vertices of $\hat{P}$ with valency 3)} .
\end{eqnarray}
By the equality (15), we see that
\[
k\leq{F-2} .
\]
Moreover, if $P$ satisfies the identity $k=F-2$, then all of the vertices of $\hat{P}$ are cusps of type $(2,2,2,2)$. 
This observation means that all of the vertices of $P$ are of type $(2,2,m)$, where $m\geq{7}$. 
Therefore, the fact that $P$ satisfies the identity $k=F-2$ contradicts to the assumption of non-compactness of $P$.
Thus, we obtain the following inequality.
\[
k\leq{F-3}.
\]

Suppose that $k=F-3$. Since $P$ is non-compact, $P$ has at least one cusp. Then if $P$ has a cusp of type $(2,3,6)$ or $(2,4,4)$ or $(3,3,3)$, $\hat{P}$ also has a cusp of  type $(2,3,6)$ or $(2,4,4)$ or $(3,3,3)$.
By the fact that a $\frac{\pi}{m}$-edge ($m\geq{3}$) is adjacent to two vertices with valency $3$,  we can see that $\hat{P}$ has at least three vertices with valency $3$. 

Therefore, by the equality (15), we obtain the following inequality.\\
\vspace{-6mm}
\begin{align*}
F-2 &= v_{2,2,2,2}+k+\dfrac{1}{2}\text{(the number of vertices of $\hat{P}$ with valency 3)} \\
       &\geq F-3+\dfrac{3}{2} = F-\dfrac{3}{2}
\end{align*}
Hence if $P$ has a cusp of type $(2,3,6)$ or $(2,4,4)$ or $(3,3,3)$, we arrive at a contradiction.\\
This implies that if $k=F-3$, $P$ has a unique cusp of type $(2,2,2,2)$ and all other vertices of $P$ are of type $(2,2,m)$  where $m\geq{7}$. \qed

\subsection{The growth rates in the case of $k=F-3$}
By Theorem 3, if $P$ satisfies that the identity $k=F-3$, we can assume that $P$ has a unique cusp of type $(2,2,2,2)$ and the other vertices are type $(2,2,m)$, where $m\geq{7}$. 
Suppose $\tilde{P}$ is the polyhedron obtained by opening all of cusps of type $(2,2,2,2)$ of $P$.
Then $\tilde{P}$ is a compact Coxeter polyhedron, therefore we can apply the following Theorem 5 to conclude that the growth rate of $P$ is a Perron number.

\begin{theo} (\cite{K}, Theorem 5.)
Suppose $\tilde{P}$ is a compact Coxeter polyhedron all of whose vertices are type of $(2,2,m)$ where $m\geq{7}$ and $P$ is the pinched Coxeter polyhedron of $\tilde{P}$ which has a unique cusp of type $(2,2,2,2)$ and the other vertices are type $(2,2,m)$, where $m\geq{7}$.\\
Then, the growth rate of $P$ is a Pisot number.
\end{theo}

Since the growth rates of hyperbolic Coxeter polyhedra with $F-3$ $\frac{\pi}{m}$-edges for $m\geq{7}$ are Perron numbers by Theorem 5, we consider growth rates of hyperbolic Coxeter polyhedra with at most $F-4$ $\frac{\pi}{m}$-edges for $m\geq{7}$ in the next subsection.

\subsection{The growth rates in the case of $k\leq{F-4}$ part 1}
In this subsection, we prove the following proposition. 
\begin{prop}
Suppose that $k\leq{F-4}$ and $P$ satisfies the following inequality 
\begin{equation}
v_{2,2,2,2}+e_{3}+e_{4}+e_{5}+e_{6}+F-8 \geq{0}.
\end{equation}
 Then the growth rate of $P$ is a Perron number.
\end{prop}

In order to prove Proposition 2, we shall use the following notation and terminology. 
\begin{defi}{(abstract polyhedron)}

An \textit{abstract polyhedron} $C$ is a simple graph $C$ on 2-dimensional sphere $S^2$ all of its vertices are 3-valent or 4-valent.
If each edge of an abstract polyhedron $C$ is labeled with $\frac{\pi}{m}$ for an integer $m\geq{2}$ or $m=\infty$, $C$ is called an \textit{abstract Coxeter polyhedron}. 
\end{defi} 

For any hyperbolic Coxeter polyhedron $P$, $\partial{P}$ is homeomorphic to $S^2$. 
This implies that the face structure of $P$ gives the abstract Coxeter polyhedron $C$. 
We call $C$ the abstract Coxeter polyhedron associated with $P$. 
Suppose that $C$ is an abstract Coxeter polyhedron and $v$ is a vertex with valency $k$ for $k=3$ or $k=4$. 
Let $c_1, \cdots, c_k$ be edges of $C$ incident to $v$ and $\frac{\pi}{a_i}$ be the label on the edge $c_i$. 
\begin{itemize}
  \item If a vertex $v$ of $C$ with valency 3 satisfies the inequality $\frac{1}{a_1}+\frac{1}{a_2}+\frac{1}{a_3}>1$, we call $v$ a \textit{spherical vertex} \textit{of type} $(a_1, a_2, a_3)$. 
  \item If a vertex $v$ of $C$ with valency 3 satisfies the equality $\frac{1}{a_1}+\frac{1}{a_2}+\frac{1}{a_3}=1$, we call $v$ a \textit{Euclidean vertex} \textit{of type} $(a_1, a_2, a_3)$. 
  \item If a vertex $v$ of $C$ with valency 3 satisfies the inequality $\frac{1}{a_1}+\frac{1}{a_2}+\frac{1}{a_3}<1$, we call $v$ a \textit{hyperbolic vertex} \textit{of type} $(a_1, a_2, a_3)$.
  \item If a vertex $v$ of $C$ with valency 4 satisfies the equality $a_1=a_2=a_3=a_4=2$, we call $v$ a \textit{Euclidean vertex of type} $(2,2,2,2)$. 
  \item A vertex $v$ of $C$ with valency 4 other than Euclidean vertex is called a \textit{hyperbolic vertex}. 
  \item $\mathcal{V}_{a_1, a_2, a_3}$ denotes the number of spherical vertices of type $(a_1, a_2, a_3)$. 
  \item $\mathcal{E}_m$ denotes the number of edges labeled by $\frac{\pi}{m}$. 
  \item $\mathcal{F}$ denotes the number of faces of $C$. 
\end{itemize}
A spherical, Euclidean or hyperbolic vertex $v$ of $C$ corresponds to a spherical, Euclidean or hyperbolic Coxeter triangle $\Delta_{a_1, a_2, a_3}$ whose interior angles are $\frac{\pi}{a_1}, \frac{\pi}{a_2}$ and $\frac{\pi}{a_3}$, respectively. 
We denote by $\Gamma_{a_1, a_2, a_3}$ the Coxeter group associated with $\Delta_{a_1, a_2, a_3}$.
The growth series of $\Gamma_{a_1, a_2, a_3}$ is denoted by $f_{a_1, a_2, a_3}(t)$.
Then the \textit{pseudo growth function $f_C(t)$} of $C$ is defined by the following identity. 
$$
\frac{1}{f_C(t^{-1})}:=1-\frac{\mathcal{F}}{[2]}+\sum _{m\geq{2}} \frac{\mathcal{E}_m}{[2, m]}-\sum _{a_1, a_2, a_3} \frac{\mathcal{V}_{a_1, a_2, a_3}}{f_{a_1, a_2, a_3}(t)}.
$$

Let $P$ be a hyperbolic Coxeter polyhedron and $C$ be the abstract Coxeter polyhedron associated with $P$. 
Then we can see that the pseudo growth function $f_{C}(t)$ of $C$ is equal to the growth function $f_{S}(t)$ of $P$. 

Suppose that $P$ has some dihedral angles $\frac{\pi}{m}$ for $m\geq{7}$ and $C$ is the abstract Coxeter polyhedron associated with $P$. 
Let $P'$ be the abstract Coxeter polyhedron obtained from $C$ by changing one of the dihedral angles  of $C$ from $\frac{\pi}{m}$ to $\frac{\pi}{6}$. 
By Andreev's theorem \cite{A}, the endpoints of a $\frac{\pi}{m}$-edge of $P$ is the vertices of type $(2, 2, m)$ for $m\geq{7}$, so that the abstract polyhedron $P'$ has at least 1 Euclidean vertices and no hyperbolic vertices. 
Then the growth function $f_S(t)$ of $P$ differs from the pseudo growth function $f_{P'}(t)$ of $P'$ in the terms related to changing the dihedral angle. 
This implies the following equality. 

\begin{figure}[htbp]
\begin{center}
 \includegraphics [width=300pt, clip]{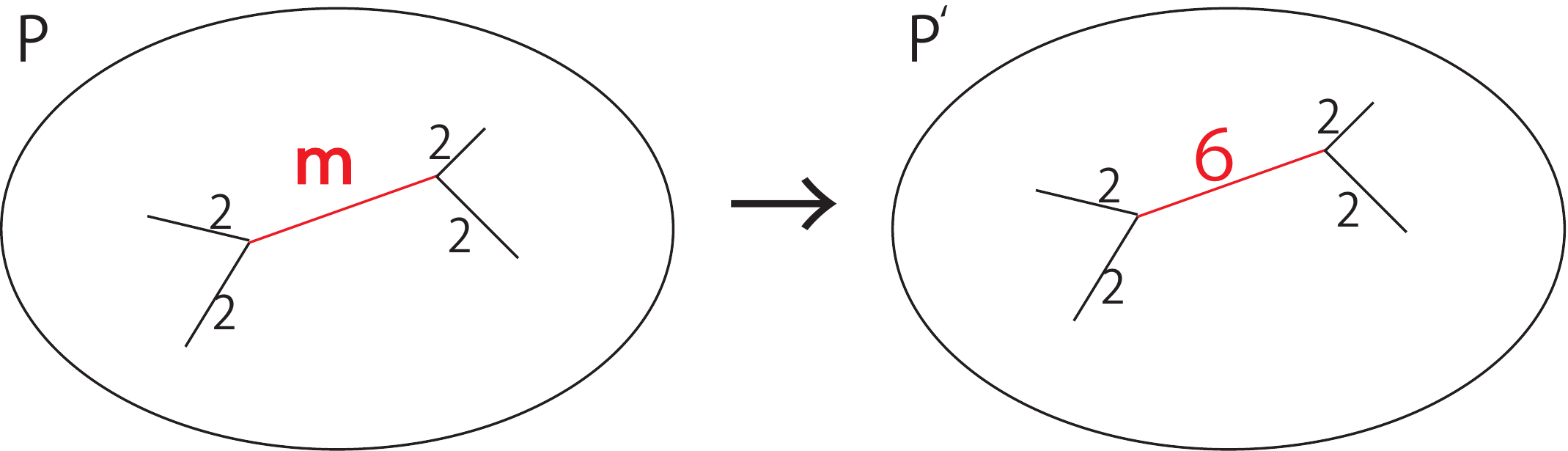}
\end{center}
\caption{}
\end{figure} 

\begin{eqnarray}
\dfrac{1}{f_{S}(t)}&=&\dfrac{1}{f_{P'}(t)}+\Bigl\{ (-\dfrac{t^6}{[2,6]}+\dfrac{2t^7}{[2,2,6]})+(\dfrac{t^m}{[2,m]}-\dfrac{2t^{m+1}}{[2,2,m]}) \Bigr\} \nonumber \\
                             &=& \dfrac{1}{f_{P'}(t)}+\dfrac{(t-1)}{[2,2,6,m]}\sum_{n=6}^{m-1}t^n
\end{eqnarray}



\proof The proof proceeds by induction on $k$. Let $P_{k}$ be a non-compact hyperbolic Coxeter polyhedron some of whose dihedral angles are $\frac{\pi}{m_{1}}, \cdots, \frac{\pi}{m_{k}}$, where $m_{i}\geq{7} \ (i=1, \cdots, k)$

\textit{Step $1$} For the case $k=1$, i.e. we consider the growth function of $P_{1}$. 
The dihedral angles of $P_{1}'$ lie in the set $\{\frac{\pi}{m} | m={2,3,4,5,6}\}$.
By the result in \cite{Y}, $\frac{1}{f_{P_{1}'}(t)}$ is written as 

\begin{eqnarray*}
\dfrac{1}{f_{P_{1}'}(t)} &=& \dfrac{(t-1)}{[4,6,10]}H_{2,3,4,5,6}(t),
\end{eqnarray*}
where $H_{2,3,4,5,6}(t)$ is the following polynomial of degree 16.
\begin{eqnarray*}
\scalebox{0.7}{
$H_{2,3,4,5,6}(t)$} &\scalebox{0.7}{=}&\scalebox{0.7}{$\left(v'_{2,3,6}+v'_{2,4,4}+v'_{3,3,3}+v'_{2,2,2,2}-1\right) t^{16}$}\\
&\scalebox{0.7}{+}&\scalebox{0.7}{$\left(F'+v'_{2,3,6}+v'_{2,4,4}+v'_{3,3,3}-4\right) t^{15}$}\\
&\scalebox{0.7}{+}&\scalebox{0.7}{$\left(\frac{1}{2} v'_{2,2,3}+\frac{1}{2}
   v'_{2,2,4}+\frac{1}{2} v'_{2,2,5}+\frac{1}{2} v'_{2,2,6}+v'_{2,3,3}+v'_{2,3,4}+v'_{2,3,5}+3 v'_{2,3,6}+3 v'_{2,4,4}+\frac{5}{2} v'_{3,3,3}+3 v'_{2,2,2,2}-4\right)
   t^{14}$}\\
   &\scalebox{0.7}{+}&\scalebox{0.7}{$\left(3 F'-\frac{1}{2} v'_{2,2,3}+\frac{1}{2} v'_{2,3,4}+\frac{1}{2} v'_{2,3,5}+\frac{5}{2} v'_{2,3,6}+2 v'_{2,4,4}+\frac{5}{2} v'_{3,3,3}-12\right)
   t^{13}$}\\
   &\scalebox{0.7}{+}&\scalebox{0.7}{$\left(\frac{3}{2} v'_{2,2,3}+v'_{2,2,4}+\frac{3}{2} v'_{2,2,5}+\frac{3}{2} v'_{2,2,6}+2 v'_{2,3,3}+\frac{5}{2} v'_{2,3,4}+3 v'_{2,3,5}+5 v'_{2,3,6}+5
   v'_{2,4,4}+\frac{9}{2} v'_{3,3,3}+5 v'_{2,2,2,2}-8\right) t^{12}$}\\
   &\scalebox{0.7}{+}&\scalebox{0.7}{$\left(5 F'-v'_{2,2,3}-\frac{1}{2} v'_{2,2,5}+v'_{2,3,4}+\frac{3}{2} v'_{2,3,5}+3 v'_{2,3,6}+3
   v'_{2,4,4}+3 v'_{3,3,3}-20\right) t^{11}$}\\
   &\scalebox{0.7}{+}&\scalebox{0.7}{$\left(2 v'_{2,2,3}+\frac{3}{2} v'_{2,2,4}+\frac{5}{2} v'_{2,2,5}+2 v'_{2,2,6}+3 v'_{2,3,3}+\frac{7}{2}
   v'_{2,3,4}+\frac{9}{2} v'_{2,3,5}+6 v'_{2,3,6}+6 v'_{2,4,4}+6 v'_{3,3,3}+6 v'_{2,2,2,2}-11\right) t^{10}$}\\
   &\scalebox{0.7}{+}&\scalebox{0.7}{$\left(6 F'-v'_{2,2,3}-v'_{2,2,5}+v'_{2,3,4}+2 v'_{2,3,5}+3
   v'_{2,3,6}+3 v'_{2,4,4}+3 v'_{3,3,3}-24\right) t^9$}\\
   &\scalebox{0.7}{+}&\scalebox{0.7}{$\left(2 v'_{2,2,3}+\frac{3}{2} v'_{2,2,4}+3 v'_{2,2,5}+2 v'_{2,2,6}+3 v'_{2,3,3}+\frac{7}{2} v'_{2,3,4}+5
   v'_{2,3,5}+6 v'_{2,3,6}+6 v'_{2,4,4}+6 v'_{3,3,3}+6 v'_{2,2,2,2}-12\right) t^8$}\\
   &\scalebox{0.7}{+}&\scalebox{0.7}{$\left(6 F'-v'_{2,2,3}-v'_{2,2,5}+v'_{2,3,4}+2 v'_{2,3,5}+3 v'_{2,3,6}+3 v'_{2,4,4}+3
   v'_{3,3,3}-24\right) t^7$}\\
   &\scalebox{0.7}{+}&\scalebox{0.7}{$\left(2 v'_{2,2,3}+\frac{3}{2} v'_{2,2,4}+\frac{5}{2} v'_{2,2,5}+2 v'_{2,2,6}+3 v'_{2,3,3}+\frac{7}{2} v'_{2,3,4}+\frac{9}{2} v'_{2,3,5}+5
   v'_{2,3,6}+5 v'_{2,4,4}+5 v'_{3,3,3}+5 v'_{2,2,2,2}-11\right) t^6$}\\
   &\scalebox{0.7}{+}&\scalebox{0.7}{$\left(5 F'-v'_{2,2,3}-\frac{1}{2} v'_{2,2,5}+v'_{2,3,4}+\frac{3}{2} v'_{2,3,5}+2 v'_{2,3,6}+2
   v'_{2,4,4}+2 v'_{3,3,3}-20\right) t^5$}\\
   &\scalebox{0.7}{+}&\scalebox{0.7}{$\left(\frac{3}{2} v'_{2,2,3}+v'_{2,2,4}+\frac{3}{2} v'_{2,2,5}+\frac{3}{2} v'_{2,2,6}+2 v'_{2,3,3}+\frac{5}{2} v'_{2,3,4}+3
   v'_{2,3,5}+3 v'_{2,3,6}+3 v'_{2,4,4}+\frac{7}{2} v'_{3,3,3}+3 v'_{2,2,2,2}-8\right) t^4$}\\
   &\scalebox{0.7}{+}&\scalebox{0.7}{$\left(3 F'-\frac{1}{2} v'_{2,2,3}+\frac{1}{2} v'_{2,3,4}+\frac{1}{2}
   v'_{2,3,5}+\frac{1}{2} v'_{2,3,6}+v'_{2,4,4}+\frac{1}{2} v'_{3,3,3}-12\right) t^3$}\\
   &\scalebox{0.7}{+}&\scalebox{0.7}{$\left(\frac{1}{2} v'_{2,2,3}+\frac{1}{2} v'_{2,2,4}+\frac{1}{2}
   v'_{2,2,5}+\frac{1}{2} v'_{2,2,6}+v'_{2,3,3}+v'_{2,3,4}+v'_{2,3,5}+v'_{2,3,6}+v'_{2,4,4}+\frac{3}{2} v'_{3,3,3}+v'_{2,2,2,2}-4\right) t^2$}\\
   &\scalebox{0.7}{+}&\scalebox{0.7}{$(F'-4) t-1$}
\end{eqnarray*} 
where $F', v'_{2,2,2,2}$ and $v'_{a_1, a_2, a_3}$ denote respectively the number of faces, vertices of type $(2,2,2,2)$ and vertices of type $(a_1, a_2, a_3)$ of $P'_1$.

Let $Q_{1}(t):=[m]H_{2,3,4,5,6}(t)+(1+2t^2+2t^4+2t^6+2t^8+t^{10})\sum_{n=6}^{m-1}t^n$. Then the degree of $Q_{1}(t)$ is $m+15$, so that we can represent $Q_{1}(t)$ as $\sum_{i=1}^{m+15}a_{i}^{(1)}t^i$. 
By the main result of \cite{Y}, all of the coefficients of $H_{2,3,4,5,6}(t)$ except its constant term are non-negative. 
Therefore we can obtain the following inequality and identities.
\begin{eqnarray*}
a_{i}^{(1)} &\geq& 0 \ (i\geq{6}) \\
a_{5}^{(1)} &=& a_{5}+a_{4}+a_{3}+a_{2}+a_{1}-1\\
a_{4}^{(1)} &=& a_{4}+a_{3}+a_{2}+a_{1}-1\\
a_{3}^{(1)} &=& a_{3}+a_{2}+a_{1}-1\\
a_{2}^{(1)} &=& a_{2}+a_{1}-1=v'_{2,2,2,2}+e'_{3}+e'_{4}+e'_{5}+e'_{6}+F'-9 \\
a_{1}^{(1)} &=& a_{1}-1 = F'-5, 
\end{eqnarray*}
where $\{a_{i}\}_{i=1}^{16}$ are the coefficients of $H_{2,3,4,5,6}(t)=\sum_{i=1}^{16}a_{i}t^{i}-1$ with respect to $P'_{1}$ and $v'_{2,2,2,2}$ and $e'_{m}$ are the number of cusps of type $(2,2,2,2)$ and $\frac{\pi}{m}$-edges of $P'_{1}$, respectively.

Since $P'_{1}$ is obtained from $P_{1}$ by changing dihedral angles from $\frac{\pi}{m}$ to $\frac{\pi}{6}$, $a_{2}^{(1)}$ can be rewritten as 
\[
a_{2}^{(1)}=v_{2,2,2,2}+e_{3}+e_{4}+e_{5}+e_{6}+F-8 .
\]

The above equality means that the coefficients of $Q_{1}(t)$ except its constant term are non-negative under the assumption of Proposition 2.

\textit{Step 2.} Suppose that the following identity holds for the growth function $f_{S_{k-1}}(t)$ of $P_{k-1}$ where $k\geq{2}$.
\[
\dfrac{1}{f_{S_{k-1}}(t)}=\dfrac{(t-1)Q_{k-1}(t)}{[2,2,6,m_{1},\cdots, m_{k-1}](1+2t^2+2t^4+2t^6+2t^8+t^{10})} , 
\] 
where $Q_{k-1}(t)$ is a polynomial of degree $m_{1}+\cdots+m_{k-1}+16-(k-1)$ and the coefficients of $Q_{k-1}(t)$ except its constant term are non-negative. By the identity (17) we deduce that the following identities.

\begin{eqnarray*}
\dfrac{1}{f_{S_{k}}(t)} &=& \dfrac{(t-1)}{[2,2,6]} \Bigl\{ \dfrac{Q_{k-1}(t)}{[m_{1},\cdots,m_{k-1}](1+2t^2+2t^4+2t^6+2t^8+t^{10})}+\dfrac{\sum_{n=6}^{m_{k}-1}t^n}{[m_{k}]} \Bigr\} \\
                               &=& \dfrac{(t-1)\{ [m_{k}]Q_{k-1}(t)+[m_{1},\cdots,m_{k-1}](1+2t^2+2t^4+2t^6+2t^8+t^{10})\sum_{n=6}^{m_{k}-1}t^n \}}{[2,2,6,m_{1},\cdots,m_{k}](1+2t^2+2t^4+2t^6+2t^8+t^{10})}
\end{eqnarray*} 

Let $Q_{k}(t):=[m_{k}]Q_{k-1}(t)+[m_{1},\cdots,m_{k-1}](1+2t^2+2t^4+2t^6+2t^8+t^{10})\sum_{n=6}^{m_{k}-1}t^n$. Then the degree of $Q_{1}(t)$ is $m_{1}+\cdots+m_{k}+16-k$, and hence we obtain the following inequalities and identities once we represent  $Q_{k}(t)$ as $\sum a_{i}^{(k)}t^i-1$.
\begin{eqnarray*}
a_{i}^{(k)} &\geq& 0 \ (i\geq{6}) .\\
a_{5}^{(k)} &=& a_{5}^{(k-1)}+a_{4}^{(k-1)}+a_{3}^{(k-1)}+a_{2}^{(k-1)}+a_{1}^{(k-1)}-1.\\
a_{4}^{(k)} &=& a_{4}^{(k-1)}+a_{3}^{(k-1)}+a_{2}^{(k-1)}+a_{1}^{(k-1)}-1.\\
a_{3}^{(k)} &=& a_{3}^{(k-1)}+a_{2}^{(k-1)}+a_{1}^{(k-1)}-1.\\
a_{2}^{(k)} &=& a_{2}^{(k-1)}+a_{1}^{(k-1)}-1.\\
a_{1}^{(k)} &=& a_{1}^{(k-1)}-1 = a_{1}^{(1)}-(k-1) . 
\end{eqnarray*}

By the result of \textit{Step 1}, 
\[
a_{1}^{(k)}=a_{1}^{(1)}-(k-1)=F-4-k
\]

Therefore the coefficients of $Q_{k}(t)$ except its constant term are non-negative if $P$ satisfies the inequality $F-4\geq{k}$, and hence the growth rate is a Perron number.
\qed

\subsection{The growth rates in the case of $k\leq{F-4}$ part 2}
In Proposition 2, the inequality (16) is sufficient to show the growth rate of $P$ is a Perron number for the case of $F-4\geq{k}$. 
Let $P$ be a non-compact hyperbolic Coxeter polyhedron with at least 7 facets. 
Then $P$ has at least 1 cusps. 
This implies the following inequality: 
\[
v_{2,2,2,2}+e_3+e_4+e_5+e_6+F-8\geq{1+7-8}=0.
\] 
For this reason, in this subsection we consider non-compact Coxeter polyhedra with 5 or 6 facets which do not satisfy the inequality (16).
The Figure 3 shows the combinatorial structures of convex polyhedra with 4, 5, 6 facets.

\begin{figure}[htbp]
\begin{center}
 \includegraphics [width=150pt, clip]{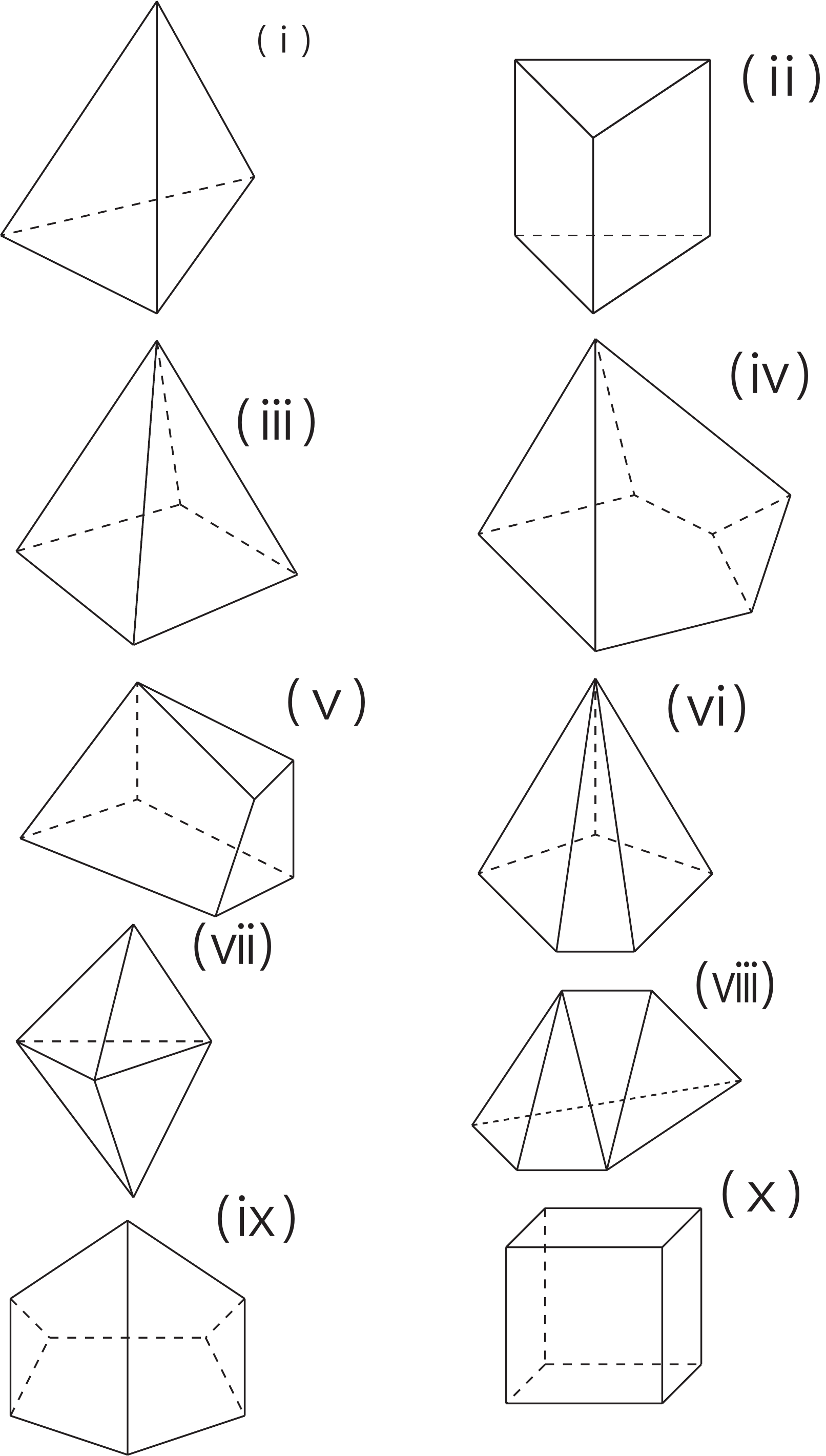}
\end{center}
\caption{}
\label{fig2}
\end{figure} 

We use Andreev's theorem to determine when the non-compact hyperbolic Coxeter polyhedron $P$ with 5 or 6 facets does not satisfy the inequality (16).

\begin{theo}(\cite{A} Andreev's theorem)
An acute-angled almost simple polyhedron of finite volume with given dihedral angles, other than a tetrahedron or a triangular prism, exists in $\mathbb{H}^{3}$ if and only if the following conditions are satisfied:

(a) if three facets meet at a vertex, then the sum of the dihedral angles between them is greater than $\pi$;

(b) if three facets meet at a cusp, then the sum of the dihedral angles between them is equal to $\pi$;

(c) if four facets meet at a vertex or a cusp, then all the dihedral angles between them equal $\frac{\pi}{2}$;

(d) if three faces are pairwise adjacent but share neither a vertex nor a cusp, then the sum of the dihedral angles between them is less than $\pi$ ;

(e) if a facet $F_{i}$ is adjacent to facets $F_{j}$ and $F_{k}$, while $F_{j}$ and $F_{k}$ are not adjacent but have a common cusp which $F_{i}$ does not share, then at least one of the angles formed by $F_{i}$ with $F_{j}$ and with $F_{k}$ is different from $\frac{\pi}{2}$;

(f) if four facets are cyclically adjacent but meet at neither a vertex nor a cusp, then at least one of the dihedral angles between them is different from $\frac{\pi}{2}$. 

\end{theo}

A hyperbolic polyhedron in $\mathbb{H}^{3}$ is almost simple if any of its finite vertices have valence 3. 
Thus any hyperbolic Coxeter polyhedron is almost simple. 

By Theorem 4 and Andreev's theorem \cite{A}, if non-compact hyperbolic Coxeter polyhedra with 5 or 6 facets have $\frac{\pi}{m}$-edges for $m\geq{7}$, then their combinatorial structures could be (ii), (iv), (v), (viii), (ix), (x).
If the combinatorial structure of $P$ is (viii), $P$ has 2 cusps of type $(2,2,2,2)$ and if the combinatorial structure is (ix) or (x), $P$ has at least one of cusps of type $(2,3,6)$ or $(2,4,4)$ or $(3,3,3)$. Hence, the inequality (16) holds for the  combinatorial structures (viii), (iv), (x).

By Theorem 4, we determine which edges can have dihedral angles $\frac{\pi}{m}$ for $m\geq{7}$. 
The possible sequences of combinatorial structures obtained by opening cusps of type $(2,2,2,2)$ are shown in Figure 4.

\begin{figure}[htbp]
   \begin{center}
      \begin{tabular}{c}
             
             \begin{minipage}{0.45\hsize}
               \begin{center}
                 \includegraphics[clip, width=5.5cm]{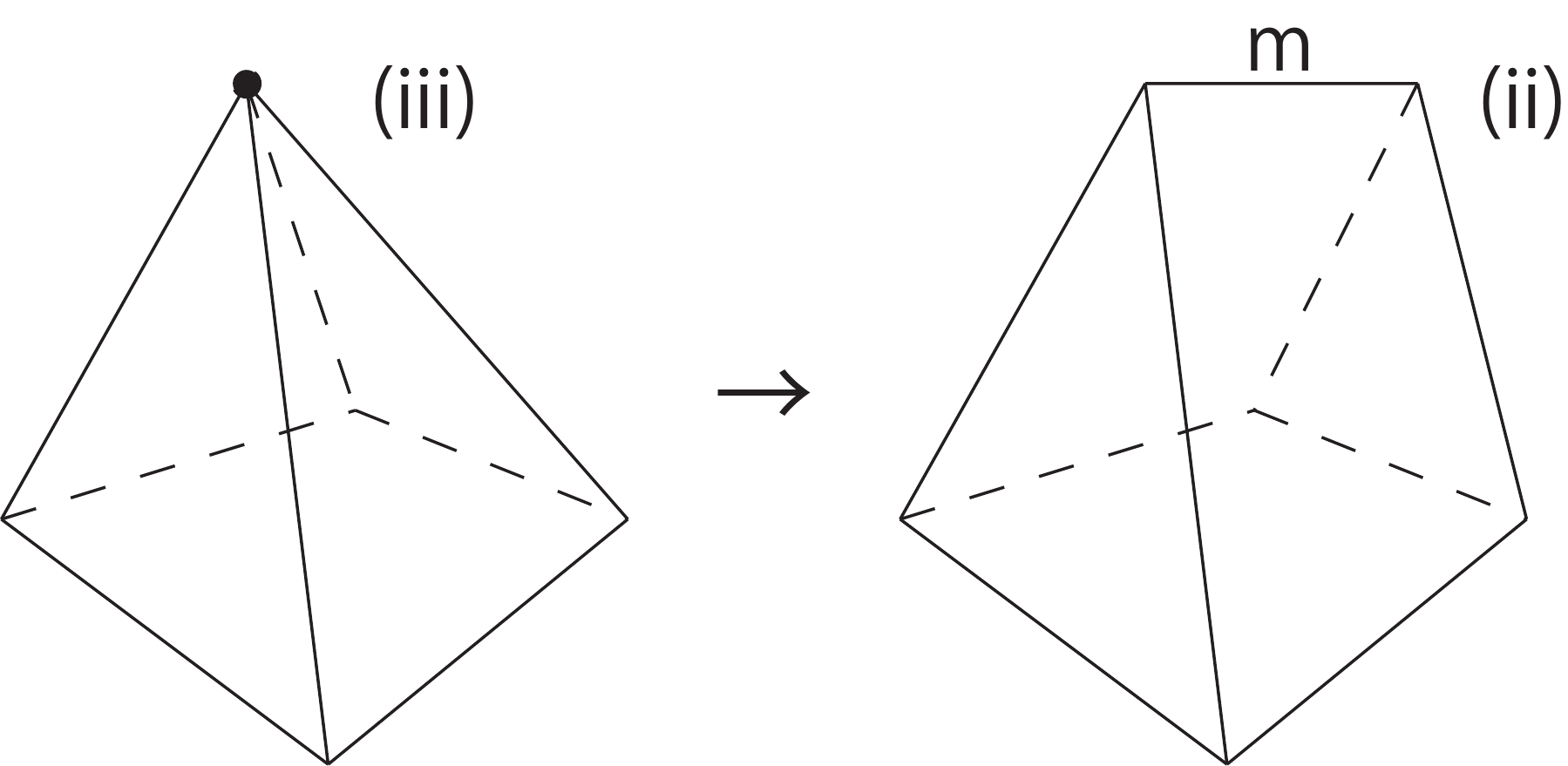}
                 \hspace{2.0cm} 
               \end{center}
             \end{minipage}
             
             \begin{minipage}{0.45\hsize}
               \begin{center}
                 \includegraphics[clip, width=5.5cm]{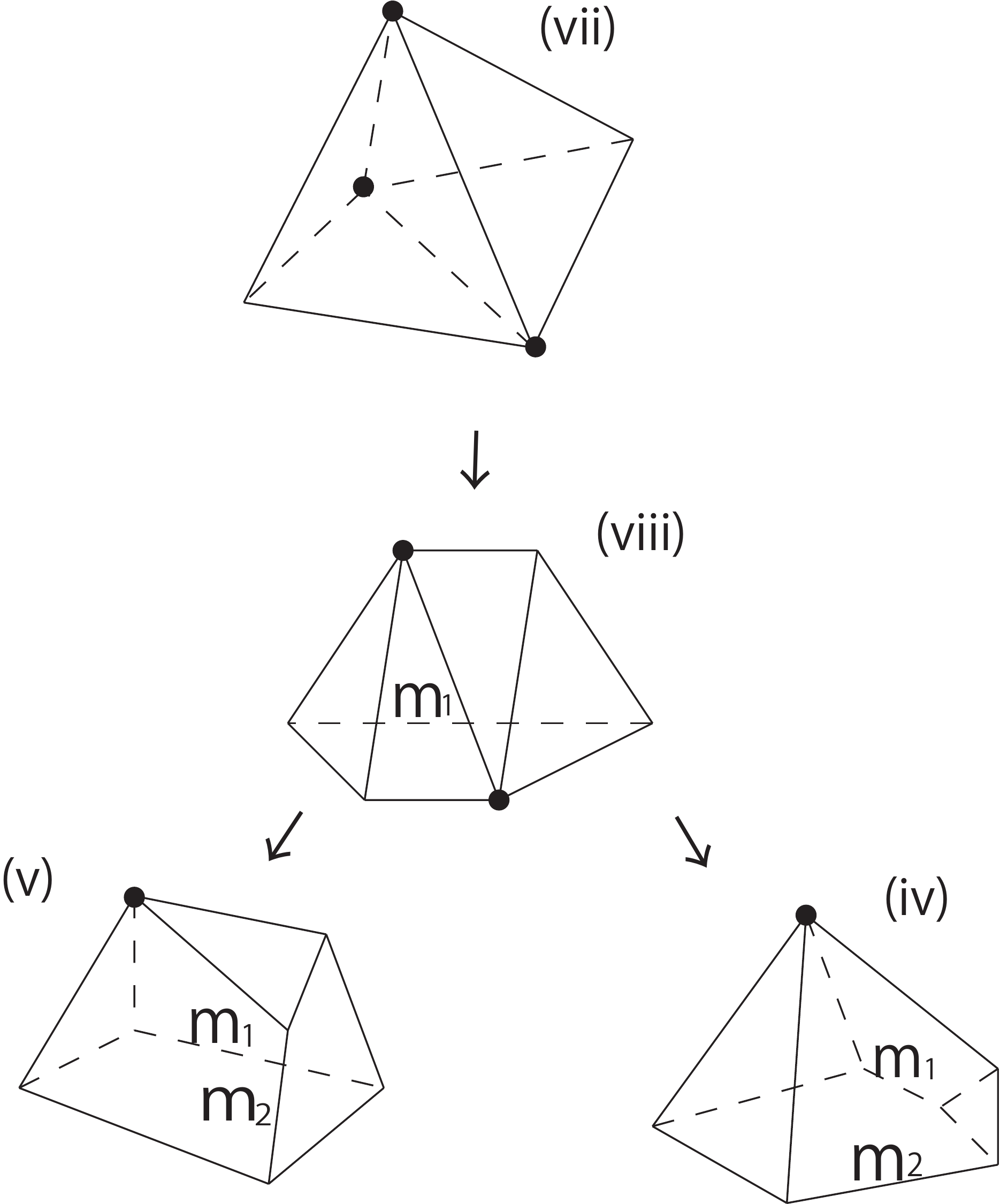}
                 \hspace{2.0cm} 
               \end{center}
             \end{minipage}

       \end{tabular}
     \end{center}
   \caption{}
   \label{fig3}
  \end{figure}

In Figure 4, labels on edges mean the dihedral angles and $m, m_{1}, m_{2}\geq{7}$.
If the inequality (16) does not hold for the case of (iv) or (v), all of the dihedral angles other than $\frac{\pi}{m_{1}}, \frac{\pi}{m_{2}}$ are $\frac{\pi}{2}$, since $v_{2,2,2,2}=1$.

\begin{prop}
Suppose that the combinatorial structure of $P$ is \textrm{(iv)} or \textrm{(v)} and $P$ does not satisfy the inequality (16). Then the growth rate of $P$ is a Perron number.
\end{prop}
\proof
By means of Steinberg's formula, we can calculate the growth function $f_{S}(t)$ of $P$ as
\begin{eqnarray*}
\dfrac{1}{f_{S}(t)} &=& 1-\dfrac{6t}{[2]}+\dfrac{9t^2}{[2,2]}+\dfrac{t^{m_{1}}}{[2,m_{1}]}+\dfrac{t^{m_{2}}}{[2,m_{2}]}-\dfrac{2t^{3}}{[2,2,2]}-\dfrac{2t^{m_{1}+1}}{[2,2,m_{1}]}-\dfrac{2t^{m_{2}+1}}{[2,2,m_{2}]}\\
                               &=& \dfrac{(t-1)\bigl\{ (2t+1)[m_{1},m_{2}]-(t+1)([m_{1}]+[m_{2}])\bigr\}}{[2,2,2,m_{1},m_{2}]}
\end{eqnarray*}
Let $Q(t):=(2t+1)[m_{1},m_{2}]-(t+1)([m_{1}]+[m_{2}])$. We may assume that $m_{1}\geq{m_{2}}$, without loss in generality.

If $m_{1}=m_{2}$,  $Q(t)$ can be rewrriten as, 
\begin{eqnarray*}
Q(t) &=& [m_{1}]\bigl\{ (2t+1)[m_{1}]-(2t+2) \bigr\} \\
       &=& [m_{1}] \bigl( 2\sum_{k=0}^{m_{1}-1}t^{k+1}+\sum_{k=0}^{m_{1}-1}t^{k}-2t-2\bigr)\\
       &=& [m_{1}](2t^{m_{1}}+3t^{m_{1}-1}+3t^{m_{1}-2}+\cdots +3t^2+t-1)
\end{eqnarray*}

If $m_{1}>m_{2}$, $Q(t)$ can be rewrriten as,
\begin{eqnarray*}
Q(t) &=& (2t+1)\Bigl\{ (t^{m_{1}-1}+\cdots+t^{m_{2}})[m_{2}]+[m_{2}]^{2} \Bigr\}-(t+1)\Bigl\{ (t^{m_{1}-1}+\cdots+t^{m_{2}})+2[m_{2}]\Bigr\}\\
       &=&(2t+1)(t^{m_{1}-1}+\cdots+t^{m_{2}})[m_{2}]-(t+1)(t^{m_{1}-1}+\cdots+t^{m_{2}})+ [m_{2}]\bigl\{ (2t+1)[m_{2}]-(2t+2) \bigr\} \\
       &=&[m_{1}](2t^{m_{2}}+3t^{m_{2}-1}+\cdots +3t^2+t)+t(t^{m_{1}-1}+\cdots+t^{m_{2}})-[m_{2}]
\end{eqnarray*} 

By the above calculation, the coefficients of $Q(t)$ except its constant term are non-negative.
 
Therefore we can apply Proposition 1 to conclude that the growth rate is a Perron number.
\qed

Coxeter triangular prisms one of whose dihedral angles is $\frac{\pi}{m}$ for $m\geq{7}$ are obtained from Coxeter pyramids by opening the cusp of type $(2,2,2,2)$. 
We use the following lemma to determine when a Coxeter triangular prism one of whose dihedral angels is $\frac{\pi}{m}$ for $m\geq{7}$ can be realized in $\mathbb{H}^{3}$.  

\begin{lem}\cite{KU2}
Coxeter pyramids in $\mathbb{H}^{3}$ can be classified by Coxeter graphs in Fig 5 up to isometry.
\end{lem}
\begin{figure}[htbp]
\begin{center}
 \includegraphics [width=250pt, clip]{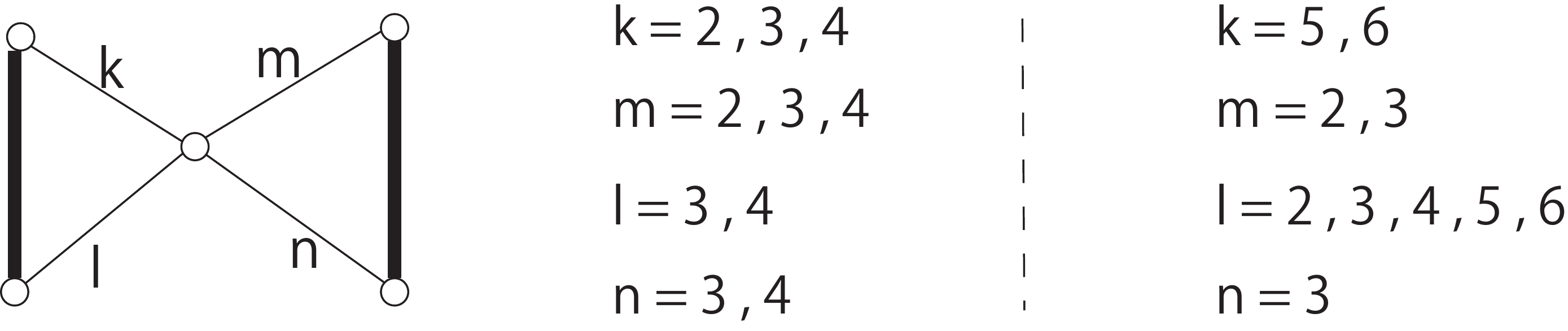}
\end{center}
\caption{}
\label{fig5}
\end{figure}

For any non-compact Coxeter triangular prism one of whose dihedral angles is $\frac{\pi}{m}$ for $m\geq{7}$, there exists a unique Coxeter pyramid which has some cusps other than the apex. 
Therefore, by Lemma 1, we can determine when Coxeter triangular prisms do not satisfy the condition (16).

\begin{cor}
Suppose that $P$ is a Coxeter triangular prism and $P$ does not satisfy the inequality (16). 
Then $P$ has the dihedral angles as in Figure 6 and the growth rate of $P$ is a Perron number.
\end{cor}
\begin{figure}[htbp]
\begin{center}
 \includegraphics [width=200pt, clip]{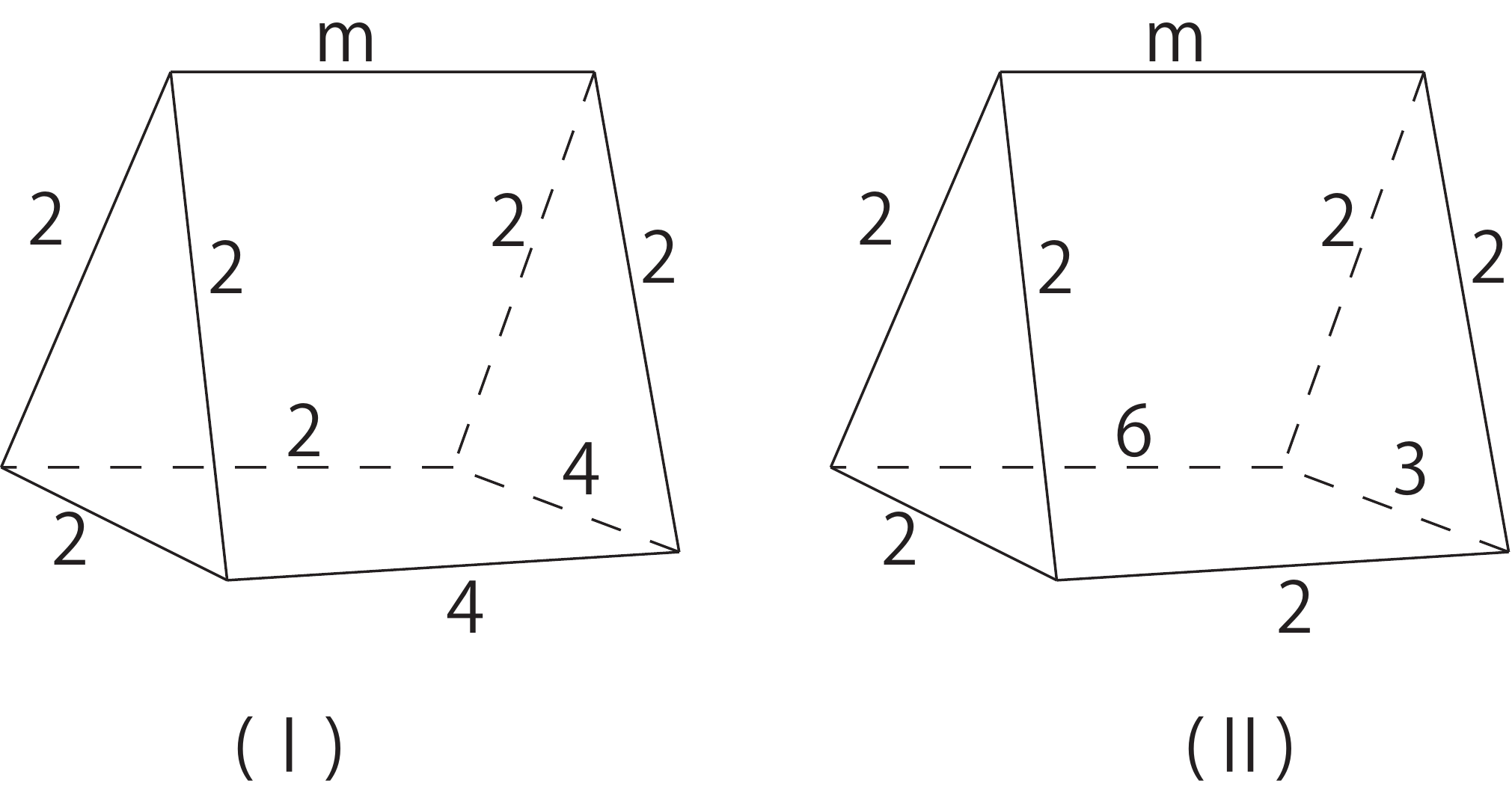}
\end{center}
\caption{}
\label{fig5}
\end{figure} 

\proof
\textit{Case.(I)} By means of Steinberg's formula, we can calculate the growth function $f_{S}(t)$ of $P$, and hence the growth function is written as,  
\begin{eqnarray*}
\dfrac{1}{f_{S}(t)} &=&\dfrac{(t-1)(2t^{m+2}+3t^{m+1}+4t^m+\cdots+4t^4+3t^3+t^2-1)}{[2,2,4,m]}  
\end{eqnarray*}

\textit{Case.(II)} The growth function is calculated in the same manner:   
\begin{eqnarray*}
\dfrac{1}{f_{S}(t)} &=&\dfrac{R(t)}{[2,2,2,3,6,m]}  
\end{eqnarray*} 
where
\begin{center}
$R(t)=2 t^{m+8}+5t^{m+7}+7 t^{m+6}+7 t^{m+5}+6t^{m+4}+5 t^{m+3}+3 t^{m+2}+t^{m+1}$\\
 \ \ \ \ \ \ \ \ \ \ \ \ \ \ \ $-t^9-4 t^8-7 t^7-8t^6-7 t^5-6 t^4-4 t^3-t^2+t+1$
\end{center}

Let us notice that $R(t)$ is divisible by [2,3] and $(t-1)$. Therefore $f_{S}(t)$ can be rewritten as, 
\begin{eqnarray*}
\dfrac{1}{f_{S}(t)} &=&\dfrac{(t-1)(2t^{m+4}+3t^{m+3}+4t^{m+2}+5t^{m+1}+6t^m+\cdots+6t^6+5t^5+3t^4+2t^3+t^2-1)}{[2,2,6,m]}  
\end{eqnarray*} 

Hence, we can apply Proposition 1 to conclude that the growth rate is a Perron number.
\qed

By combining Theorem 5, Proposition 2-3 and Corollary 1, we obtain Theorem A. 
The classical result by Parry in \cite{Pa}, together with the result in \cite{Y} and Theorem A imply Theorem B. 


\section{Acknowledgement}
The author thanks Professor Yohei Komori for helpful comments when calculating the growth functions of non-compact Coxeter polyhedra with some dihedral angles $\frac{\pi}{m}$ for $m\geq{7}$. 



\begin{thebibliography}{}
\bibitem{A}
E. M. Andreev, On convex polyhedra of finite volume in Lobachevskij space, Mat. Sb., Nov. Ser. {\bf 83} (1970), 256-260. English transl.: Math. USSR, Sb. {\bf 12} (1971), 255- 259.
\bibitem{dlH1}
P. de la Harpe, Groupes de Coxeter infinis non affines, Exposition. Math {\bf5} (1987),  91--96.
\bibitem{Hu}
J. E. Humphreys, {\em Reflection groups and Coxeter groups}, Cambridge Studies in Advanced Mathematics, {\bf 29}, Cambridge Univ. Press, Cambridge, 1990.
\bibitem{Ke}
R. Kellerhals, Cofinite hyperbolic Coxeter groups, minimal growth rate and Pisot numbers, Algebr. Geom. Topol. {\bf 13} (2013), 1001-1025.
\bibitem{KK}
R. Kellerhals and A. Kolpakov, The minimal growth rate of cocompact Coxeter groups in hyperbolic 3-space, Canad. J. Math. {\bf 66}(2014), 354-372.
\bibitem{KP}
R. Kellerhals and G. Perren, On the growth of cocompact hyperbolic Coxeter groups, 
European J. Combin. {\bf 32} (2011), no. 8, 1299-1316.
\bibitem{K}
A. Kolpakov, Deformation of finite-volume hyperbolic Coxeter polyhedra, limiting growth rates and Pisot numbers, 
European J. Combin. {\bf33}, (2012), 1709-1724.
\bibitem{KU}
Y. Komori and Y. Umemoto, On the growth of hyperbolic 3-dimensional generalized simplex reflection groups,
Proc. Japan Acad. Ser. A Math. Sci. Volume {\bf 88}, Number 4 (2012), 62--65.
\bibitem{KU2}
Y. Komori and Y. Umemoto, On 3-dimensional hyperbolic Coxeter pyramids,
arXiv: 1503.00583
\bibitem{KY}
Y. Komori and T. Yukita, On the growth rate of ideal Coxeter groups in hyperbolic 3-space, 
Proc. Japan Acad. Ser. A Math. Sci. Volume {\bf 91}, Number 10 (2015), 155--159.
\bibitem{KN}
J. Nonaka and R.Kellerhals, The growth rates of ideal Coxeter polyhedra in hyperbolic 3-space,
To appear in Tokyo Journal of Mathematics.
\bibitem{Pa}
W. Parry, Growth series of Coxeter groups and Salem numbers, 
J. Algebra {\bf154}, (1993), 406-415
\bibitem{So}
L. Solomon, The orders of the finite Chevalley groups, J. Algebra {\bf 3}(1966), 376--393.
\bibitem{St}
R. Steinberg, {\em Endomorphisms of linear algebraic groups}, Memoirs of the American Mathematical Society, No. 80, Amer. Math. Soc., Providence, RI, 1968.
\bibitem{Y}
T.Yukita, On the growth rates of cofinite 3-dimensional hyperbolic Coxeter groups whose dihedral angles are of the form $\frac{\pi}{m}$ for $m=2,3,4,5,6$, To appear in RIMS K\^oky\^uroku Bessatsu, arXiv:1603.04592v4.

\end{thebibliography}
\end{document}